\magnification\magstep1
\hsize 15.5truecm
\scrollmode

 at 17.28truept
 at 14.4truept
 at 12truept
 at 10.95truept
 at 17.28truept
 at 10truept

\def\title#1{%
\vskip0pt plus.3\vsize\penalty-100%
\vskip0pt plus-.3\vsize\bigskip\vskip\parskip%
\bigbreak\bigbreak\centerline{\bf #1}\bigskip%
}

\def\chapter#1#2{\vfill\eject
\centerline{\bf Chapter #1}
\vskip 6truept%
\centerline{\bf #2}%
\vskip 2 true cm}

\def\section#1#2{%
\def\\{#2}%
\vskip0pt plus.3\vsize\penalty-100%
\vskip0pt plus-.3\vsize\bigskip\vskip\parskip%
\par\noindent{\bf #1\hskip 6truept%
\ifx\empty\\{\relax}\else{\bf #2\smallskip}\fi}}

\def\subsection#1#2{%
\def\\{#2}%
\vskip0pt plus.3\vsize\penalty-20%
\vskip0pt plus-.3\vsize\medskip\vskip\parskip%
\def\TEST{#1}%
\noindent{\ifx\TEST\empty\relax\else\bf #1\hskip 6truept\fi%
\ifx\empty\\{\relax}\else{#2\smallskip}\fi}}

\def\proclaim#1{\medbreak\begingroup\noindent{\bf #1.---}\enspace\sl}

\def\endproclaim{\endgroup\par\medbreak}

\def\lbreak{\hfil\break}


\def\comfig#1#2\par{
\medskip
\centerline{\hbox{\hsize=10cm\eightpoint\baselineskip=10pt
\vbox{\noindent #1}}}\par\centerline{ Figure #2}}

\def\figcom#1#2\par{
\medskip
\centerline
{Figure #1}
\par\centerline{\hbox{\hsize=10cm\eightpoint\baselineskip=10pt
\vbox{\noindent #2}}}}
\def\bull{~\vrule height .9ex width .8ex depth -.1ex}


\def\lbreak{\hfill\break}

\def\comfig#1#2\par{
\medskip
\centerline{\hbox{\hsize=10cm\eightpoint\baselineskip=10pt
\vbox{\noindent{\sl  #1}}}}\par\centerline{{\bf Figure #2}}}

\def\figcom#1#2\par{
\medskip
\centerline
{{\bf Figure #1}}
\par\centerline{\hbox{\hsize=10cm\eightpoint\baselineskip=10pt
\vbox{\noindent{\sl  #2}}}}}

\def\em{\sl}

\def\\{\hfill\break}

\def\lbreak{\hfill\break}

\def\bull{~\vrule height .9ex width .8ex depth -.1ex}

\def\a{\alpha}
\def\b{\beta}
\def\g{\gamma}

\def\d{\delta}

\def\La{\Lambda}

\def\CC{{\bf C}}

\def\cD{{\cal D}}

\def\cM{{\cal M}}


\def\sqr#1#2{{\vcenter{\hrule height.#2pt%
\hbox{\vrule width.#2pt height#1pt\kern#1pt%
\vrule width.#2pt}%
\hrule height.#2pt}}}

\newfam\gotfam
\font\twlgot=eufm10 at 12pt
\font\tengot=eufm10

\font\sevengot=eufm7
\textfont\gotfam=\twlgot
\scriptfont\gotfam=\tengot 
\scriptscriptfont\gotfam=\sevengot


\input amssym.def
\input amssym.tex

\overfullrule=0pt
\input psfig
\input psbox
\psfordvips

\def\deg{{\rm deg}}\def\Ad{{\rm
Ad}} \def\deg{{\rm deg}} 
  \def\gg{{\bf g}}  

\centerline{\bf Dynamical 
$r$-matrices for Hitchin's systems on Schottky curves}
\bigskip
\centerline{B. Enriquez}
\medskip
{\bf Abstract.} {\em
We express Hitchin's systems on curves in Schottky parametrization, and 
construct dynamical $r$-matrices attached to them.
}
\medskip

\section{}{Introduction.} This paper is devoted to the study of
Hitchin's integrable systems. These systems are defined on the cotangent
space to the moduli space of holomorphic $G$-bundles on a Riemann
surface. One source of the recent interest in these systems is that they
can be viewed as a classical limit of the $\cD$-modules proposed by
Beilinson and Drinfeld in their approach to the geometric Langlands
correspondence ([BD]). 

A classical approach to integrable systems is the $r$-matrix approach;
this viewpoint usually leads to their quantization. In a previous
work ([ER]), we constructed such $r$-matrices in the case of Hitchin
systems in genus one. These $r$-matrices are dynamical; it means that
they depend on a (Poisson commutative) set of phase space variables. Due
to this dependence, they satisfy a dynamical generalization of the classical
Yang-Baxter equation (DYBE). In his work [F], G. Felder derived the same
$r$-matrices from considerations related to the
Knizhnik-Zamolodchikov-Bernard equation. Moreover, he constructed a
quantum version of the generalized Yang-Baxter equation and quantum
groups attached to them. 

In this paper, we formulate Hitchin's integrable systems on a higher
genus Riemann surface, uniformized \`a la Schottky, and point out a
connection with the Garnier systems in the case of Mumford curves. (Such
a formulation can also be found in [O].)
We then derive a
dynamical $r$-matrix for these systems. This gives an alternative proof of
their integrability. We also show that they satisfy an analogue of the
DYBE, but it seems difficult to formulate correctly a quantum
counterpart of this equation. 

Let us also note here that an adelic analogue of the $r$-matrix presented
here has been found by G. Felder in [Fe]. It would be interesting to
relate both constructions. We also hope that the version presented here
of the DYBE may help to understand the correct "higher genus" version of
the formalism of P. Etingof and A. Varchenko ([EV]). 

We express our thanks to J. Harnad, Y. Kosmann-Schwarzbach, 
V. Rubtsov and A. Reyman, with whom
we had discussions related to the subject of this work, and to the
former for telling us about the work of Garnier.  
 
\section{}{1. Hitchin's systems with level structures.}

We recall this notion from [Be], [Ma].
Let $X$ be a compact complex curve, $D=\sum_{i}[P_{i}]$ a divisor on $X$
(all points $P_{i}$ are assumed different). Let $G$ be a complex reductive
group, ${\bf g}$ be its Lie algebra.
Let $\cM_{G,D}(X)$ the moduli space of pairs $(P,j_{i})$ of a
principal $G$-bundle over $X$ and of trivialisations of $P$ above each
point $P_{i}$. Then the cotangent space to $\cM_{G,D}(X)$ at $(P,j_{i})$
is identified with $H^{0}(X,{\bf g}_{P}\otimes\Omega^{1}(-D))$, where
$\gg_{P}=P\times_{G}\gg$ ($\gg$ being viewed as the adjoint representation of
$G$), and $\Omega^{1}$ being the canonical bundle of $X$. 
We have a natural mapping ([Hi])
$T^{*}\cM_{G,D}(X)\to\oplus_{i=1}^{r}H^{0}(X,\Omega^{1}(-D)^{\otimes
d_{i}})$, $r$ being the rank of $G$ and $d_{i}$ the exponents of $G$
($d_{i}=\deg P_{i}$, for $(P_{i})_{1\le i\le r}$ 
a basis of the invariant polynomials on
$\gg$), obtained by
applying $P_{i}$ on $H^{0}(X,g_{P}\otimes\Omega^{1}(-D))$. Functions on 
$T^{*}\cM_{G,D}(X)$ obtained from this mapping are in involution,
$T^{*}\cM_{G,D}(X)$ being endowed with its natural symplectic structure. 

\section{}{2. Schottky uniformisation.}

Let $l\ge 1$ be an integer, and let on $\CC P^{1}$, $\Gamma_{i}$ and
$\Gamma'_{i}$ ($i=1,\cdots, l$) be $2l$ circles, bounding $2l$ disjoint open
discs $D_{i}$ and $D'_{i}$. Let us give ourselves $l$ elements
$\gamma_{i}$ of $SL(2,\CC)$ mapping $\Gamma_{i}$ to $\Gamma'_{i}$; and let
$X$ be the Riemann surface 
$$X=\CC
P^{1}-\bigcup_{i=1}^{l}D_{i}\cup D'_{i}/(x_{i}\sim\gamma_{i}(x_{i}),
x_{i}\in\Gamma_{i}).$$ 

An open subset of the moduli space $\cM_{G}(X)$ of principal
$G$-bundles over $X$ can be identified with $G^{l}/G$ (where $G$ acts on
each factor of $G^{l}$ by the adjoint action), by associating to
the class of $(g_{i})_{1\le i\le l}$ the class of the bundle 
$
P_{(g_{i})}=\big[\CC
P^{1}-\bigcup_{i=1}^{l}D_{i}\cup D'_{i}\big]\times G /((x_{i},g)
\sim(\gamma_{i}(x_{i}),
g_{i}g), x_{i}\in\Gamma_{i}, g\in G).
$

Assume for simplicity, that $\infty$ belongs to $\CC
P^{1}-\bigcup_{i=1}^{l}D_{i}\cup D'_{i}$. Then an open subset of 
$\cM_{G,[\infty]}(X)$ can be identified with $G^{l}$:
adjoin to the bundle defined above, the identity mapping from the fibre at
$\infty$ to $G$. 

Let us identify now the cotangent bundle to $G$ with $G\times \gg^{*}$ via
left invariant one-forms, and accordingly $T^{*}G^{l}$ with
$G^{l}\times \gg^{*l}$. Elements of $H^{0}(X,\gg_{P}\otimes\Omega^{1}
(-[\infty]))$ are
given by twisted holomorphic one-forms, with poles at $\infty$, that is
by one-forms over $\CC P^{1}-\La$ ($\La$ is the limit set of the free group
$\Gamma$ generated by the $\gamma_{i}$) $\omega(z)dz$, such that 
$$
\omega(\gamma z)\g'(z)dz=\Ad(g_{\g})\omega(z)dz, \leqno{(1)}
$$
for all $\g\in\Gamma$ and with only poles at $\Gamma\infty$. Recall the
definition of the multiplier $q_{\g}$ of an hyperbolic element $\g$ of
$SL(2,\CC)$: the transformation defined by $\gamma$ verifies
${{\g(z)-a_{\g}}\over{\g(z)-b_{\g}}}
=q_{\g}{{z-a_{\g}}\over{z-b_{\g}}}$, for certain $a_{\g}$ and
$b_{\g}$ in $\CC P^{1}$ and $q_{\g}\in\CC^{\times}$, $|q_{\g}|< 1$. 
We have: 

\proclaim{Lemma 1}
Over the open subset of $G^{l}$ defined by
$$\sum_{i=1}^{l}q_{\g_{i}}(\|\Ad 
g_{i}\|+ \|\Ad
g_{i}^{-1}\|)\lbreak < 1,
$$
($\|\  \|$ an algebra norm in End $g$), 
 the twisted holomorphic one-form corresponding
to $(g_{i},\xi_{i})\in G^{l}\times \gg^{*l}$ is given by the Poincar\'e
series (cf. also [B]) 
$$
\xi(g_{i},\xi_{i})(z)dz=
\sum_{\g\in\Gamma, 1\le i\le l}{{\Ad(g_{\g}^{-1})\xi_{i}}\over{\gamma(z)-
\gamma_{i}^{-1}(\infty)}}\gamma'(z)dz. \leqno{(2)}
$$
\endproclaim
\medskip
\noindent 
\bf Proof. \rm Let us prove the 
convergence of (2) under the present hypothesis (see also
[Bu], [Fo]). 
Let $\g=\g_{i_{1}}^{\epsilon_{1}}...\g_{i_{p}}^{\epsilon_{p}}$,
$i_{k}=1,...,l$, $\epsilon_{i_{k}}=\pm 1$. 
The norm of $\g'(z)$ is then bounded by $Cq_{i_{1}}...q_{i_{p}}$ ($C$ a
fixed constant). On
the other hand, the norm of $\Ad (g_{\g})\xi_{i}$ is estimated by 
$$
\|\Ad(g_{i_{1}}^{\epsilon_{1}})\|...\|\Ad(g_{i_{p}}^{\epsilon_{p}})\|
\|\xi_{i}\|;$$
it follows that the contribution to (2) of elements of $\Gamma$, of
length equal to $p$, is bounded by 
$$
C\Big(\sum_{i=1}^{l}q_{i}(\|\Ad g_{i}\|+\|\Ad g_{i}^{-1}\|)\Big)^{p}
(\sum_{i=1}^{l}\|\xi_{i}\|);
$$
and the sum of all these terms is bounded by 
$$
(1-\sum_{i=1}^{l}q_{i}(\|\Ad g_{i}\|+\|\Ad g_{i}^{-1}\|))^{-1}
(\sum_{i=1}^{l}\|\xi_{i}\|).
$$
We  prove similarly that the sum (2) is an analytic function of $z$.

Let $(X_{i})_{1\le i\le l}$ be $l$ elements of the Lie
algebra of $G$, and
let us view them as infinitesimal left-invariant translations. The
pairing between $(\xi_{i})$ and $(X_{i})$ is simply given by 
$
\langle (\xi_{i}), (X_{i})\rangle=\sum_{i=1}^{l}\langle \xi_{i},X_{i}
\rangle.
$
On the other hand, the pairing between $\xi(z)dz$ and $(X_{i})$ is given
by 
$
\sum_{i=1}^{l}{1\over{2i\pi}}
\int_{\Gamma_{i}}\langle X_{i},\xi(z)dz\rangle.
$
Recall that the poles of the one-form
${{\gamma'(z)dz}\over{\gamma(z)-\gamma_{i}^{-1}(\infty)}}$ are located at
$\gamma^{-1}(\infty)$ and $\gamma^{-1}\gamma_{i}^{-1}(\infty)$. The only
possibility for them to be on two different sides of $\Gamma_{i}$ is
$\gamma=e$. The contribution of the corresponding term is,
by the residues formula, $\langle \xi_{i}, X_{i} \rangle$. \bull

For the rest of the paper, we will work in the open subset defined in
prop. 1. 

\noindent \it Remarks. \rm 

1) The present formalism can be adapted to the
situation of Mumford curves, where the $q_{\g_{i}}$ are considered as formal
variables. The base ring is then $R=\CC[[q_{\g_{1}}, \cdots,\break
q_{\g_{l}}]]$. If the $g_{\g_{i}}$ belong to $G(R)$, and the $\xi_{i}$
to $\gg\otimes R$, 
the series (2) converges without restrictions; it then has to
be interpreted as a formal series of the type $dz\sum_{1\le i\le l, n\ge 1}
{{\a_{i,n}}\over{(z-a_{\g_{i}})^{n}}}+{{\b_{i,n}}\over{(z-b_{\g_{i}})^{n}}}
$, $\a_{i,n}$, $\b_{i,n}$ in $\gg\otimes R$ and with valuation 
$\ge n$, subject to
the conditions $\sum_{i=1}^{l}a_{\g_{i}}\a_{i,1}+b_{\g_{i}}\b_{i,1}=0$ 
of regularity at $\infty$. Note that Lax operators of this type
(with bounded orders of poles) already appeared in the work of Garnier
[Ga].  

2) We cannot use  ${{\g'(z)dz}\over{\g(z)-\g_{i}^{-1}(\infty})}=
{{dz}\over{z-(\g_{i}\g)^{-1}(\infty)}}-{{dz}\over{z-\g^{-1}(\infty)}}$
to regroup terms in (2) and obtain for expression of the
series (2), $\sum_{\g\in\Gamma}\Ad g_{\g}^{-1}(\sum_{i=1}^{l}\Ad 
g_{\g_{i}}\xi_{i}-\xi_{i}){{dz}\over{z-\g^{-1}(\infty)}}$, since the
last expression does not converge; but we can deduce from it the
variation of (2) under replacement of $\infty$ by a point $z_{0}$ close 
to it. The proper replacement of (2) is then  
$
\sum_{\g\in\Gamma, 1\le i\le l}\Ad(g_{\g}^{-1})\xi_{i}
({ {\gamma'(z)}\over{\gamma(z)-
\gamma_{i}^{-1}(z_{0})} }-{{\gamma'(z)}\over{\gamma(z)-z_{0}}})dz,
$
whose derivative w.r.t. $z_{0}$ is $\sum_{\g\in\Gamma}\Ad g_{\g}^{-1}
(\sum_{i=1}^{l}\big(\Ad 
g_{\g_{i}}\xi_{i}-\xi_{i})\big)
{({\g^{-1})'(z_{0})dz}\over{(z-\g^{-1}(z_{0}))^{2}}}$.
In particular, this variation is zero under the condition 
$\sum_{i=1}^{l}(\Ad 
g_{\g_{i}}\xi_{i}-\xi_{i})=0$, which is also the condition for (2) to be
regular at $\infty$ (and also the condition that the image of
$(g_{i},\xi_{i})$ by the
moment map associated to the adjoint action of $G$ on $G^{l}$ is zero). 
 
\medskip

\section{}{3. Dynamical $r$-matrices.}

Consider the ``$r$-matrices'' 
$$\eqalign{
r(z,w)dz=\sum_{\g\in\Gamma}{{\Ad
g_{\g}^{(2)}P}\over{\g(z)-w}}\g'(z)dz, \quad
s(z,w)dw=\sum_{\g\in\Gamma}{{\Ad
g_{\g}^{(1)}P}\over{z-\g(w)}} & \g'(w)dw\cr & =-r(w,z)^{(21)}dw,
}\leqno{(3)}
$$
where $P$ is the split Casimir element of $\gg\otimes \gg$ (differential
elements $dz$ and $dw$ will be considered to commute together). The
proof of convergence of analyticity of these series is similar to that
of the series (2). 

We wish to prove: 

\proclaim{Proposition} The Poisson brackets of operators $\xi$ are given
by 
 $$\eqalign{
\{\xi(g_{i},\xi_{i})(z)^{(1)}dz,\xi(g_{i},\xi_{i})(w)^{(2)}dw\}=
[r(z,w)dz, & \xi(g_{i},\xi_{i})(w)^{(2)}dw] \cr & +
[s(z,w)dw,\xi(g_{i},\xi_{i})(z)^{(1)}dz].}\leqno{(4)}
$$
\endproclaim

For this we first compare the transformation properties of both sides of
(4), under the action of $\g_{i}$ on $z$. 
Call the l.h.s. and r.h.s. of (4)
respectively $A(z,w)dzdw$ and $B(z,w)dzdw$, and set 
$C(z,w)dzdw=A(z,w)dzdw-B(z,w)dzdw$. We have: 
\proclaim{Lemma 2} $C(z,w)dzdw$ satisfies
$$
C(\g z,w)\g'(z)dzdw=\Ad g_{\g}^{(1)}C(z,w)dzdw, \leqno{(5)}
$$
$$
C(z,\g w)\g'(w)dzdw=\Ad g_{\g}^{(2)}C(z,w)dzdw, \leqno{(6)}
$$
$$
C(w,z)dzdw=-C(z,w)^{(21)}dzdw; \leqno{(7)}
$$
moreover, $C(z,w)dzdw$ has no poles on $(\CC P^{1}-\La)^{2}$. \endproclaim
\medskip
\noindent \bf Proof. \rm 
We have 
$$\eqalign{
A(\g_{i}z,w)\g_{i}'(z)dzdw =\Ad  & g_{i}^{(1)}A(z,w)dzdw
\cr & +[\{g_{i}^{(1)},\xi(w)^{(2)}dw\}g_{i}^{-1(1)},\Ad
g_{i}^{(1)}\xi(z)^{(1)}dz], 
}$$
and since
$$\eqalign{
r(\g_{i}z,w)\g_{i}'(z)dz=\Ad g_{i}^{(1)}r(z,w)dz, \ \  
s(\g_{i}z,w)dw= & \Ad g_{i}^{(1)}s(z,w)dw \cr & +\sum_{\g\in\Gamma}
{{\Ad g_{\g_{i}\g}^{(1)}P  \g'(w)dw  }\over{\g(w)-\g_{i}^{-1}(\infty)}}
}$$
(the second identity is obtained using $d_{w}\ln(\g_{i}\g(w)-\g_{i}(z))
-d_{w}\ln(\g(w)-z)=d_{w}\ln(\g(w)-\g_{i}^{-1}(\infty))$), we obtain 
$$
\eqalign{
B(\g_{i}z,w)\g_{i}'(z)dzdw= & \Ad g_{i}^{(1)}B(z,w)dzdw \cr
& +[\sum_{\g\in\Gamma}{{\Ad g_{\g_{i}\g}^{(1)}P}\over
{\g(w)-\g_{i}^{-1}(\infty)}}\g'(w)dw, 
\Ad g_{i}^{(1)}\xi(z)^{(1)}dz]. 
}$$
We then compute 
$\{g_{i}^{(1)},\xi(w)^{(2)}dw\}g_{i}^{-1(1)}$, and find using 
$\{g_{i}^{(1)},\xi_{j}^{(2)}\}=\delta_{ij}g_{i}^{(1)}P$
that it is equal to $\sum_{\g\in\Gamma}{{\Ad g_{\g_{i}\g}^{(1)}P}\over
{\g(w)-\g_{i}^{-1}(\infty)}}\g'(w)dw$. This proves (5).
(7) is clear, and together with (5) it implies
(6). 

Let us turn to the statement about poles. Let us fix
$w\notin\Gamma\infty$. Then $A(z,w)dzdw$ has poles when
$z\in\Gamma\infty$. We have for $z$ near $\infty$,
$\xi(z)dz=\sum_{i=1}^{l}(\Ad g_{i}\xi_{i}-\xi_{i}){{dz_{\infty}}\over
z_{\infty}}+$reg., 
$z_{\infty}=1/z$ is a local coordinate near $\infty$. 
Then we have 
$$\{\sum_{i=1}^{l}(\Ad
(g_{i})\xi_{i}-\xi_{i})^{(1)},\xi(w)^{(2)}dw\}
=-[P,\xi(w)^{(2)}dw],
$$
so that $A(z,w)dzdw=-[P,\xi(w)^{(2)}dw]{{dz_{\infty}}\over
z_{\infty}}+$reg. The expansion of the first part of
$B(z,w)dzdw$ is the same since
$r(z,w)dz
=-P{{dz_{\infty}}\over{z_{\infty}}}$+reg. near $z=\infty$. The second
part of $B(z,w)dzdw$ is regular since $\xi(z)dz$ has a pole of order
one, whereas $s(z,w)dw$ tends to zero as $z\to\infty$. So $C(z,w)dzdw$
has no poles for $z\to\infty$; it has no poles either for $w\to\infty$
because of (7), and these are all the poles of $A(z,w)dzdw$. The poles 
of $B(z,w)dzdw$ are the same, with the possible addition of $z\in\Gamma
w$. Because of (5), we can restrict ourselves 
to the study of the pole at $z=w$;
but this pole does not occur because of the $g$-invariance of $P$. 
\bull

We now show: 
\proclaim{Lemma 3} $C(z,w)dzdw=0$, so that (4) is valid. \endproclaim

\noindent \bf Proof. \rm Due to our assumptions on $(g_{i})$, we know
that any twisted (by $(g_{i})$) holomorphic one-form on $X$, with
possible pole at $[\infty]$ can be
written in the form of a Poincar\'e series (since these series converge,
and the dimension of the vector space they form is equal to the
dimension of $H^{0}(X,P_{(g_{i})}\times_{G}\gg(-[\infty]))$ 
as it can be computed
using the Riemann-Roch formula). It follows that there exist elements
$C_{ij}$ of $\gg\otimes\gg$, such that $C(z,w)dzdw=\sum_{1\le i,j\le
l, \g,\d\in\Gamma}
{{\g'(z)dz}\over{\g(z)-\g_{i}^{-1}(\infty)}}
{{\d'(w)dw}\over{\d(w)-\g_{j}^{-1}(\infty)}}\Ad g_{\g}^{(1)}\Ad
d_{\g}^{(2)}C_{ij} 
$. $C_{ij}$ can be computed by 
$$
C_{ij}=\int_{\Gamma_{i}\times\Gamma_{j}}C(z,w)dzdw. 
$$
Now, we have
$\int_{\Gamma_{i}\times\Gamma_{j}}A(z,w)dzdw=\{\xi_{i}^{(1)},
\xi_{j}^{(2)}\}=\d_{ij}[P,\xi_{i}^{(1)}]$, and
$$
\eqalign{
\int_{\Gamma_{i}\times\Gamma_{j}}B(z,w)dzdw & =
\int_{\Gamma_{i}\times\Gamma_{j}}
[r(z,w)dz,\xi(g_{i},\xi_{i})(w)^{(2)}dw]+ \cr
 & \int_{\Gamma_{i}\times\Gamma_{j}}
[s(z,w)dw,\xi(g_{i},\xi_{i})(z)^{(1)}dz]
,}
$$
 for $i\ne j$. In 
$\int_{\Gamma_{i}\times\Gamma_{j}}
[r(z,w)dz,\xi(g_{i},\xi_{i})(w)^{(2)}dw]$, 
we integrate first w.r.t. $z$; expanding $r(z,w)dz$
according to (3), we have to integrate the
one-form $d_{z}\ln(\g(z)-w)$; it has poles at $\g^{-1}(w)$ and
$\g^{-1}(\infty)$; $w$ being on $\Gamma_{j}$, these two points are
always on one and the same side of $\Gamma_{i}$, so that this term does
not contribute to the integral. Exchanging the roles of $z$ and $w$, we
find that $\int_{\Gamma_{i}\times\Gamma_{j}}
[s(z,w)dw,\xi(g_{i},\xi_{i})(z)^{(1)}dz]$ is also equal to
zero. Finally, $C_{ij}=0$ for $i\ne j$. 

For $i=j$, we consider a deformation $\Gamma_{i}^{\epsilon}$
of $\Gamma_{i}$, encircling $\Gamma_{i}$ and within the domain $\CC
P^{1}-\bigcup_{i=1}^{l}D_{i}\cup D'_{i}$. We have still, 
$C_{ii}=\int_{\Gamma_{i}\times\Gamma_{i}^{\epsilon}}C(z,w)dzdw$, and 
$$\int_{\Gamma_{i}\times\Gamma_{i}^{\epsilon}}A(z,w)dzdw\break
=[P,\xi_{i}^{(1)}].
$$
Now, 
$$\eqalign{
\int_{\Gamma_{i}\times\Gamma_{i}^{\epsilon}}B(z,w)dzdw & =
\int_{\Gamma_{i}\times\Gamma_{i}^{\epsilon}}
[r(z,w)dz,\xi(g_{i},\xi_{i})(w)^{(2)}dw] \cr 
& +
\int_{\Gamma_{i}\times\Gamma_{i}^{\epsilon}}
[s(z,w)dw,\xi(g_{i},\xi_{i})(z)^{(1)}dz].
} \leqno{(8)}
$$ 
Repeating the reasoning above, and due to the relative configurations of 
$\Gamma_{i}$ and $\Gamma_{i}^{\epsilon}$, the first term of the
r.h.s. of (8) is again zero. For the second term, we have
$$
\int_{\Gamma_{i}\times\Gamma_{i}^{\epsilon}}
[s(z,w)dw,\xi(g_{i},\xi_{i})(z)^{(1)}dz]=
\sum_{\g\in\Gamma} \int_{\Gamma_{i}\times\Gamma_{i}^{\epsilon}}
[{{\Ad g_{\g}^{(1)}P}\over{z-\g(w)}}\gamma'(w)dw,\xi(z)^{(1)}dz]
$$
and the only non zero contribution is from the term with $\g=1$; it
gives
$$\int_{\Gamma_{i}}dz[\int_{\Gamma_{i}^{\epsilon}}{{Pdw}\over{z-w}},
\xi(z)^{(1)}]=\int_{\Gamma_{i}}dz[P,\xi(z)^{(1)}]=[P,\xi_{i}^{(1)}].
$$
So
$C_{ii}=0$. 
\bull
\medskip
\noindent \it Remarks. \rm 

1. It is interesting to give a purely algebraic
meaning to the $r$-matrices of (3). Setting
$\rho_{\g_{i}}=\sum_{\g\in\Gamma}{
{\Ad (g_{\g}^{-1}g_{\g_{i}}^{-1})^{(2)}P
}\over{\g(w)}-\g_{i}^{-1}(\infty)}\g'(w)dw\in \gg\otimes
H^{0}(X,\Omega^{1}\otimes \gg_{P}(-[\infty]))
$ and defining $\rho_{\g}$, for 
$\g\in\Gamma$ by the
cocycle condition $\rho_{\g\g'}=\Ad g_{\g'}^{(1)}\rho_{\g}+\rho_{\g'}$,
we have
constructed a $1$-cocycle $\rho\in H^{1}(X,g_{P}\otimes H^{0}(X,
\Omega^{1}\otimes g_{P}(-[\infty])))
=H^{1}(X,\gg_{P})\otimes H^{0}(X,
\Omega^{1}\otimes \gg_{P}(-[\infty]
))$. (By Serre duality we have a natural element in 
$H^{1}(X,\gg_{P})\otimes H^{0}(X,
\Omega^{1}\otimes \gg_{P})$ but $\rho$ seems a little different.)

$\rho$ then serves to define an affine spaces bundle over the vector
bundle $H^{0}(X,\Omega^{1}\otimes \gg_{P}(-[\infty]))\otimes \gg_{P}$ over
$X$, and also over the bundle with fiber at $x\in X$, 
$H^{0}(X,\Omega^{1}\otimes \gg_{P}(-[\infty]-[x]))\otimes \gg_{P}$; 
$r$ can then be viewed as a
section of the twist of this last bundle. Probably, there is a natural
twist of
$(\Omega^{1}\otimes \gg_{P})\boxtimes \gg_{P}(-[\infty])$ over $X\times X$,
for which $r$ could be considered as a section. 

2. The result of [BV] states the local existence of an $r$-matrix for
general integrable systems; but it is easy to see that the $r$-matrices
constructed in the situation of Hitchin system, using the methods of
this work, would depend on $(\xi_{i})$ (and in fact not be defined for
$\xi_{i}=0$).  

\section{}{4. Jacobi identity}

Writing that (4) satisfies the Jacobi identity, we find the dynamical
Yang-Baxter equation: 
$$\eqalign{ & 
[r^{32}dz_{3},r^{21}dz_{2}]  +[r^{23}dz_{2},r^{31}dz_{3}]+[r^{31}dz_{3},
r^{21}dz_{2}]
+\cr
 & \sum_{\g\in\Gamma,1\le i\le l}{{\Ad
g_{\g}^{-1}(e^{\a(3)})\g'(z_{3})dz_{3}}
\over{\g(z_{3})-\g_{i}^{-1}(\infty)}}\partial_{e_{\a}^{(i)}}(r^{21}dz_{2})
\cr & +
\sum_{\g\in\Gamma,1\le i\le l}{{\Ad
g_{\g}^{-1}(e^{\a(2)})\g'(z_{2})dz_{2}}
\over{\g(z_{2})-\g_{i}^{-1}(\infty)}}\partial_{e_{\a}^{(i)}}(r^{31}dz_{3})
=0,
\cr}\leqno{(9)}$$
where $P=e^{\a}\otimes e_{\a}$, $\partial_{x^{(i)}}$ is the vector
field on $G^{l}$, given by the left translation by $x\in g$
on the $i$-th factor and zero on the other factors.

\vskip 1truecm
\noindent
{\bf References}
\bigskip
\item{[BV]} O. Babelon, C.-M. Viallet, {\sl Hamiltonian structures and
Lax equations,} Phys. Lett. B, 237 (1990), 411-6 . 

\item{[Be]} A. Beauville, CIRM lectures, Luminy (1995).

\item{[BD]} A. Beilinson, V. Drinfeld, CIRM lectures, Luminy (1995).

\item{[B]} D. Bernard, {\sl On the WZW model on Riemann surfaces,}
Nucl. Phys. B, 309 (1988), 145-174.
\item{[Bu]} W. Burnside, Proc. London Math. Soc. 23 (1891), 49. 

\item{[ER]} B. Enriquez, V. Rubtsov, {\sl  Hitchin systems, higher
Gaudin operators and $r$-matrices,} preprint alg-geom/9503010. 

\item{[EV]} P. Etingof, A. Varchenko, Geometry and classification of
solutions of the classical dynamical Yang-Baxter equation,
q-alg/9703040. 

\item{[Fo]} L.R. Ford, {\sl Automorphic functions,} Chelsea, 1951. 

\item{[Fe]} G. Felder, The KZB equations on Riemann surfaces, Les
Houches lectures, hep-th/9609153. 

\item{[Ga]} R. Garnier, Sur une classe de syst\`emes diff\'erentiels
ab\'eliens d\'eduits de la th\'eorie des \'equations lin\'eaires,
Rend. Circ. Mat. Palermo 43 (1919), 155-91.

\item{[H]} N. Hitchin, {\sl Stable bundles and integrable systems,} Duke
Math. Jour., 54:1 (1987), 91-114.

\item{[Ma]} E. Markman, {\sl Spectral curves and integrable systems,}
Comp. Math. 93 (1994), 255-290.

\item{[O]} M.A. Olshanetsky, {\sl Generalized Hitchin systems and
Knizhnik-Zamolodchikov-Bernard equation
on elliptic curves,}  Proceedings of the
6th Conference on Mathematical Physics (Moscow, June, 1995);
SISSA 125/95/EP; ITEP-TH 3/95; preprint hep-th/9510143.
\medskip
\medskip\medskip
\section{}{}

B.E.: Centre de Math\'{e}matiques, URA 169 
du CNRS, Ecole Polytechnique, 91128 Palaiseau, France

\bye